\documentclass[12pt]{article}
\usepackage{amsmath}
\usepackage{geometry} % see geometry.pdf on how to lay o ut the page. There's lots.
\title{An equation satisfied by all non-trivial zeros $\rho$ of the Riemann zeta function $\zeta$}
\author{Roupam Ghosh}
\date{\today} % delete this line to display the current date

\begin{document}

\maketitle 
\begin{center}
\textit{Abstract:}
\\
We show that if $\rho$ is a non-trivial zero of the Riemann zeta function $\zeta$
then $$2^\rho + \frac{1}{\rho - 1} + \frac{1}{2} = \rho \int_{1}^{\infty} \left\{ t + \frac{1}{2} \right\} t^{-\rho-1} dt$$
\\where, $\{ x \}$ is the fractional part of $x$.
\end{center}
\section{Introduction}
We study the Riemann zeta and Dirichlet eta function using the theory of fractional parts. By theory of fractional parts, I mean to imply a study of integrals of the form
$$\frac{1}{r} \int_{p}^{q} \left \{ f(x) \right \} dx$$
where the notation $\{ x \}$ is the fractional part of $x$. For complex numbers, $z = x + iy$ we define the fractional part as $\{ z \} = \{ x \} + i \{ y \}$. It is not hard to grasp that for all real numbers $0 < \{x\}< 1$ and for all complex numbers $0 < | \{ z \}| < \sqrt{2}$
\section{Formula's for $\zeta(s)$ and $\eta(s)$}
We already have,
\begin{equation}
\zeta(s) = \frac{s}{s-1} -  s\int_{1}^{\infty} \left \{t\right\} t^{-s-1} dt
\end{equation}
which, is known to be valid for all $\Re(s) > 0$.  [1]
\\\\
We derive a similar integral that gives the Dirichlet eta function $\eta(s)$ for all $\Re(s) > 0$
\subsection{An expression for $\eta(s)$ using fractional part}
\textbf{Theorem:}
\textit{For the Dirichlet eta function defined by for all $\Re(s) > 0$}
$\eta(s) = \sum_{k = 1}^{\infty} \frac{1}{(2k - 1)^s} - \frac{1}{(2k)^s}$
\textit{we get an equivalent expression in the form for all $\Re(s) > 0$
where $ \kappa(t) = \left\{ t/2 \right\} + 1/2 -  \left\{ t/2 + 1/2 \right\} $, and the expression being given by 
\begin{equation}
\eta(s) = s\int_{1}^{\infty} \kappa(t) t^{-s-1} dt
\end{equation}
}
\textbf{Proof:}
A simple simplification of the integral shall prove this case. We have for all $t \in [1,\infty]$ $\kappa(t) = 0$ or $1$. It is not hard to see that $\kappa(t) = 0$ whenever $t \in [2k,2k+1)$ and $1$ whenever $t \in [2k-1,2k)$ for all positive integers $k$. Hence, we can write the integral as
$$s\int_{1}^{\infty} \kappa(t) t^{-s-1} dt = \sum_{k=1}^{\infty} \int_{2k-1}^{2k} s t^{-s-1} dt$$
Giving us the sum
$$\sum_{k = 1}^{\infty} \frac{1}{(2k - 1)^s} - \frac{1}{(2k)^s}$$
which is nothing but $\eta(s)$. Since we already know that this sum converges for $\Re(s) > 0$, we get our result.
\section{The equation for non-trivial zeros}
\textbf{Theorem:} \textit{The non-trivial zeros $\rho$ of the Riemann zeta function $( 0 < \Re(\rho) < 1 )$ satisfy the equation
$$
2^\rho + \frac{1}{\rho - 1} + \frac{1}{2} = \rho \int_{1}^{\infty} \left\{ t + \frac{1}{2} \right\} t^{-\rho-1} dt
$$
}
\textbf{Proof:}
We know that if $\zeta(\rho) = 0$ then so is $\eta(\rho) = 0$ for $0 < \Re(\rho) < 1$. [2] From equation (1) we get
\begin{equation}
\frac{\rho}{\rho - 1} = \rho \int_{1}^{\infty} \left \{t\right\} t^{- \rho -1} dt
\end{equation}
Now from equation (2) substituting $2t$ in place of $t$ and simplifying for $\eta(\rho) = 0$ gives us
\begin{equation}
\rho \int_{1/2}^{1} \kappa(2t) t^{-\rho-1} dt
+ \rho \int_{1}^{\infty} \frac{1}{2}t^{-\rho-1} dt
 + \rho \int_{1}^{\infty} \left\{ t \right\} t^{-\rho-1} dt
 = \rho \int_{1}^{\infty} \left\{ t + \frac{1}{2} \right\} t^{-\rho-1} dt
\end{equation}
Whenever $t \in [1/2,1)$ we have $\kappa(2t) = 1$. Substituting (3) in equation (4) and evaluating gives us
\begin{equation}
2^\rho + \frac{1}{\rho - 1} + \frac{1}{2} = \rho \int_{1}^{\infty} \left\{ t + \frac{1}{2} \right\} t^{-\rho-1} dt
\end{equation}
which is our desired result.
\\\\
Equation (5) is the main equation in this paper. All non-trivial zeros $\rho$ of the Riemann zeta function satisify this equation. \section{Acknowledgements}
As always, my thanks goes to my mom, dad, and sister for making my life wonderful, and also, to my friend Craig Feinstein.


\begin{thebibliography}{99} 
\bibitem{1}Peter Borwein, Stephen Choi, Brendan Rooney, Andrea Weirathmueller \\
\textit{"The Riemann Hypothesis: A Resource for the Afficionado and Virtuoso Alike"}. (pg. 12) 
\bibitem{2}Peter Borwein, Stephen Choi, Brendan Rooney, Andrea Weirathmueller \\
\textit{"The Riemann Hypothesis: A Resource for the Afficionado and Virtuoso Alike"}. (pg. 49) 
\end{thebibliography}
\end{document}